\documentclass{amsart}
\usepackage{amsmath, amsthm, amscd, amsfonts, amssymb, graphicx, color}
\usepackage[bookmarksnumbered, colorlinks, plainpages]{hyperref}

\newtheorem{theorem}{Theorem}[section]

\newtheorem{corollary}[theorem]{Corollary}
\theoremstyle{definition}
\newtheorem{definition}[theorem]{Definition}
\newtheorem{example}[theorem]{Example}

\theoremstyle{remark}

\numberwithin{equation}{section}
\newcommand{\ff}{if and only if }

\begin{document}

\title[Arens regularity and weak amenability]{Arens regularity of module actions\\ and weak amenability of Banach algebras}
\author[]{KAZEM HAGHNEJAD AZAR }

\dedicatory{}



\begin{abstract}For Banach left and right module actions,  we will establish  the relationships between  topological centers of module actions with some result in the weak amenability of Banach algebras.
\end{abstract}

\maketitle

\section{Preliminaries and Introduction}
\noindent As is well-known [1], the second dual $A^{**}$ of Banach algebra $A$ endowed with the either Arens multiplications is a Banach algebra. The constructions of the two Arens multiplications in $A^{**}$ lead us to definition of topological centers for $A^{**}$ with respect to both Arens multiplications. The topological centers of Banach algebras, module actions and applications of them  were introduced and discussed in [5,  9, 14,  17]. In this paper, we will study some properties of topological centers of module actions and we offer some relations between factorization properties and topological properties of Banach $A-bimodule$ $X$ with topological centers properties of module actions.
For Banach algebra $A$, Dales,  Rodrigues-Palacios and  Velasco in [7] have been studied the weak amenability of $A$, when its second dual is weakly amenable.   Mohamadzadih and Vishki in [17] have given simple solution to this problem with some other results, and  Eshaghi Gordji and Filali in [10] have been studied  this problem with some new results. We study this problem in the new way with some new results.\\
Now we introduce some notations and definitions that we used
throughout this paper as follows.\\
 Let $X$ be a normed space and $X^*,\,X^{**}$ be the first and the second dual of $X,$ respectively. We use $(X^*, \sigma(X^*,X))$ and $(X, \sigma(X,X^*))$ to denote the spaces $X^*$ and $X$ with respect to the weak$^*$ and weak topology, respectively.

Let $\mathcal{A}$ be  a Banach algebra. For $a\in \mathcal{A}$
 and $a^\prime\in \mathcal{A}^*$, we denote by $a^\prime a$
 and $a a^\prime,$ respectively, the functionals in $\mathcal{A}^*$ defined by $\langle a^\prime a,b\rangle =\langle a^\prime,ab\rangle$ and $\langle a a^\prime,b\rangle =\langle a^\prime,ba\rangle,$ for all $b\in \mathcal{A}$.
   The Banach algebra $\mathcal{A}$ is embedded in its second dual via the identification
 $\langle a,a^\prime\rangle =\langle a^\prime,a\rangle $ for every $a\in
\mathcal{A}$ and $a^\prime\in
\mathcal{A}^*$. A bounded net $(e_{\alpha})_{{\alpha}\in I}$ in $\mathcal{A}$ is called a \emph{bounded left
approximate identity} (BLAI) [resp. \emph{bounded right
approximate identity} (BRAI)] if,
 for each $a\in \mathcal{A}$,   $e_{\alpha}a\rightarrow a$ [resp. $ae_{\alpha}\rightarrow a$]. Moreover, $(e_{\alpha})$ is called a (two sided) bounded approximate identity (BAI), if for every $a\in \mathcal{A},$ the conditions $e_{\alpha}a\rightarrow a$ and $ae_{\alpha}\rightarrow a$ both hold.

  Let $X,Y$ and $Z$ be normed spaces and $m:X\times Y\rightarrow Z$ be a bounded bilinear map. Richard Arens in \cite{1}, offers two natural extensions $m^{***}$ and $m^{t***t}$ of $m$ from $X^{**}\times Y^{**}$ into $Z^{**}$ as follows. First we give $m^{***}$ in three steps:
\begin{enumerate}
\item[(1)] $m^*:Z^*\times X\rightarrow Y^*$,~~~~~given by~~~$\langle m^*(z^\prime,x),y\rangle =\langle z^\prime, m(x,y)\rangle $, where $x\in X$, $y\in Y$ and $z^\prime\in Z^*;$
 \item[(2)] $m^{**}:Y^{**}\times Z^{*}\rightarrow X^*$,~~given by $\langle m^{**}(y^{\prime\prime},z^\prime),x\rangle =\langle y^{\prime\prime},m^*(z^\prime,x)\rangle,$ where $x\in X$, $y^{\prime\prime}\in Y^{**}$ and $z^\prime\in Z^*;$
 \item[(3)] $m^{***}:X^{**}\times Y^{**}\rightarrow Z^{**}$,~ given by~ ~ ~$\langle m^{***}(x^{\prime\prime},y^{\prime\prime}),z^\prime\rangle $  $=\langle x^{\prime\prime},m^{**}(y^{\prime\prime},z^\prime)\rangle $\\ ~where ~$x^{\prime\prime}\in X^{**}$, $y^{\prime\prime}\in Y^{**}$ and $z^\prime\in Z^*$.
\end{enumerate}
The map $m^{***}$ is the unique extension of $m$ such that the mapping $x^{\prime\prime}\mapsto m^{***}(x^{\prime\prime},y^{\prime\prime})$ from $X^{**}$ into $Z^{**}$ is weak$^*$ to weak$^*$ continuous for every $y^{\prime\prime}\in Y^{**}$, and the mapping $y^{\prime\prime}\mapsto m^{***}(x,y^{\prime\prime})$ is weak$^*$ to weak$^*$ continuous from $Y^{**}$ into $Z^{**},$ for every $x\in X$. Then the \emph{first topological center} of $m$ is defined as
$$Z_1(m)=\{x^{\prime\prime}\in X^{**}:~~y^{\prime\prime}\rightarrow m^{***}(x^{\prime\prime},y^{\prime\prime}) \text{ is weak}^* \text{ to weak}^* \text{ continuous}\}.$$

Let now $m^t:Y\times X\rightarrow Z$ be the \emph{transpose} of $m$ defined by $m^t(y,x)=m(x,y),$ for all $x\in X$ and $y\in Y$. Then $m^t$ is a continuous bilinear map and so it may be extended, as above, to the mapping $m^{t***}:Y^{**}\times X^{**}\rightarrow Z^{**}$.
 The mapping $m^{t***t}:X^{**}\times Y^{**}\rightarrow Z^{**}$ in general is not equal to $m^{***}$, see \cite{1}; if it happens, then we say that $m$ is \emph{Arens regular}. The mapping $y^{\prime\prime}\mapsto m^{t***t}(x^{\prime\prime},y^{\prime\prime})$ is weak$^*$ to weak$^*$ continuous for every $y^{\prime\prime}\in Y^{**}$, but the mapping $x^{\prime\prime}\mapsto m^{t***t}(x^{\prime\prime},y^{\prime\prime})$ from $X^{**}$ into $Z^{**}$ is not in general  weak$^*$ to weak$^*$ continuous for every $y^{\prime\prime}\in Y^{**}$. So we define the \emph{second topological center} of $m$ as
$$Z_2(m)=\{y^{\prime\prime}\in Y^{**}:~~x^{\prime\prime}\rightarrow m^{t***t}(x^{\prime\prime},y^{\prime\prime}) \text{ is weak}^*\text{ to weak}^* \text{ continuous}\}.$$
It is clear that $m$ is Arens regular if and only if $Z_1(m)=X^{**},$ or equivalently $Z_2(m)=Y^{**}$. Arens regularity of $m$ is also equivalent to the following condition
$$\lim_i\lim_j\langle z^\prime,m(x_i,y_j)\rangle =\lim_j\lim_i\langle z^\prime,m(x_i,y_j)\rangle ,$$
whenever both limits exist for all bounded sequences $(x_i)_i\subseteq X$, $(y_i)_i\subseteq Y$ and $z^\prime\in Z^*$, see \cite{6,15,19}. On the other hand, the map $m$ is called \emph{left strongly Arens irregular} if $Z_1(m)=X,$ and \emph{right strongly Arens irregular} if $Z_2(m)=Y$.

Let $\mathcal{A}$ be a normed algebra with the algebra multiplication $\pi.$ Then the regularity of $\mathcal{A}$ is defined to be the regularity of the map $\pi,$ when is considered as a bilinear map. Let $a^{\prime\prime}$ and $b^{\prime\prime}$ be elements of $\mathcal{A}^{**}.$ By the Goldstin's theorem \cite[p.425]{8}, there are nets $(a_{\alpha})_{\alpha}$ and $(b_{\beta})_{\beta}$ in $\mathcal{A}$ such that $a_\alpha \rightarrow a^{\prime\prime}$ and  $b_\beta \rightarrow b^{\prime\prime},$ both in $(\mathcal{A}^{**}, \sigma(\mathcal{A}^{**},\mathcal{A}^{*})).$ So it is easy to see that for all $a^\prime\in \mathcal{A}^*$,
$$\lim_{\alpha}\lim_{\beta}\langle a^\prime,\pi(a_{\alpha},b_{\beta})\rangle =\langle a^{\prime\prime}\cdot b^{\prime\prime},a^\prime\rangle $$ and
$$\lim_{\beta}\lim_{\alpha}\langle a^\prime,\pi(a_{\alpha},b_{\beta})\rangle =\langle a^{\prime\prime}ob^{\prime\prime},a^\prime\rangle ,$$
where $a^{\prime\prime}\cdot b^{\prime\prime}$ and $a^{\prime\prime}ob^{\prime\prime}$ are the first and second Arens products of $\mathcal{A}^{**}$, respectively, see \cite{6,15,19}.

For a Banach algebra $\mathcal{A}$ with the operation $\pi:\mathcal{A}\times \mathcal{A}\rightarrow \mathcal{A}$ , we shall also simplify our notations. The first (left) Arens product of $a^{\prime\prime},b^{\prime\prime}\in \mathcal{A}^{**}$ shall be simply defined by the three steps:
 $$\langle a^\prime a,b\rangle =\langle a^\prime ,ab\rangle ,$$
  $$\langle a^{\prime\prime} a^\prime,a\rangle =\langle a^{\prime\prime}, a^\prime a\rangle ,$$
  $$\langle a^{\prime\prime}\cdot b^{\prime\prime},a^\prime\rangle =\langle a^{\prime\prime},b^{\prime\prime}a^\prime\rangle .$$
 for every $a,b\in \mathcal{A}$ and $a^\prime\in \mathcal{A}^*$. Similarly, the second (right) Arens product of $a^{\prime\prime},b^{\prime\prime}\in \mathcal{A}^{**}$ shall be defined by:
 $$\langle a oa^\prime ,b\rangle =\langle a^\prime ,ba\rangle ,$$
  $$\langle a^\prime oa^{\prime\prime} ,a\rangle =\langle a^{\prime\prime},a oa^\prime \rangle ,$$
  $$\langle a^{\prime\prime}ob^{\prime\prime},a^\prime\rangle =\langle b^{\prime\prime},a^\prime ob^{\prime\prime}\rangle .$$
  for all $a,b\in A$ and $a^\prime\in \mathcal{A}^*$. We remark that $\mathcal{A}^{**}$ with respect to either Arens products is a Banach algebra.

We find the usual first and second topological center of $\mathcal{A}^{**}$, which are
  $$Z_1(\mathcal{A}^{**})=Z_1(\pi)=\{a^{\prime\prime}\in \mathcal{A}^{**}: b^{\prime\prime}\rightarrow a^{\prime\prime}\,\cdot\,b^{\prime\prime} \text{ is weak}^* \text{ to weak}^* \text{ continuous}\},$$
   $$Z_2(\mathcal{A}^{**})=Z_2(\pi)=\{a^{\prime\prime}\in \mathcal{A}^{**}: a^{\prime\prime}\rightarrow a^{\prime\prime}ob^{\prime\prime} \text{ is weak}^* \text{ to weak}^* \text{ continuous}\}.$$

 An element $e^{\prime\prime}$ of $\mathcal{A}^{**}$ is said to be a \emph{mixed unit,} if $e^{\prime\prime}$ is a
right unit for the first Arens multiplication and a left unit for
the second Arens multiplication. That is, $e^{\prime\prime}$ is a mixed unit if
and only if,
for each $a^{\prime\prime}\in \mathcal{A}^{**}$, $a^{\prime\prime}\cdot e^{\prime\prime}=e^{\prime\prime}o a^{\prime\prime}=a^{\prime\prime}$. By
\cite[p.146]{4}, an element $e^{\prime\prime}$ of $\mathcal{A}^{**}$  is  mixed
      unit if and only if it is a weak$^*$ cluster point of some BAI $(e_\alpha)_{\alpha \in I}$  in
      $\mathcal{A}$.

Let now $X$ be a Banach $\mathcal{A}$-bimodule, and let
$$\pi_\ell:~\mathcal{A}\times X\rightarrow X\quad\text{and}\quad\pi_r:~X\times \mathcal{A}\rightarrow X$$
be the right and left module actions of $\mathcal{A}$ on $X$. Then $X^*$ is a left Banach $\mathcal{A}$-module and a right Banach $\mathcal{A}$-module with respect the module actions $\pi_r^{t*t}$ and $\pi_\ell^*,$ respectively. Also the second dual
$X^{**}$ is a Banach $\mathcal{A}^{**}$-bimodule with the module actions
$$\pi_\ell^{***}:~\mathcal{A}^{**}\times X^{**}\rightarrow X^{**}\quad\text{and}\quad\pi_r^{***}:~X^{**}\times \mathcal{A}^{**}\rightarrow X^{**},$$
where $\mathcal{A}^{**}$ is considered as a Banach algebra with respect to the first Arens product.
Similarly, $X^{**}$ is a Banach $\mathcal{A}^{**}$-bimodule with the module actions
$$\pi_\ell^{t***t}:~\mathcal{A}^{**}\times X^{**}\rightarrow X^{**}\quad\text{and}\quad\pi_r^{t***t}:~X^{**}\times \mathcal{A}^{**}\rightarrow X^{**},$$
where $\mathcal{A}^{**}$ is considered as a Banach algebra with respect to the second Arens product.

We say that a Banach right $\mathcal{A}$-bimodule $X$ with right module action $\pi_r,$ \emph{factors on the left} with respect to $\mathcal{A},$ if $\pi_r$ is onto, i.e., $\pi_r(X,\mathcal{A})=X.$ Similarly, if $X$ is left Banach $\mathcal{A}$-module and $\pi_\ell(\mathcal{A},X)=X,$ we say that $X$ \emph{factors on the right}, with respect to $\mathcal{A}.$ In particular, we say that $\mathcal{A}^*$ factors on the left (right), if $\mathcal{A}^*=\mathcal{A}^*\mathcal{A}$ ($\mathcal{A}^*=\mathcal{A}\mathcal{A}^*$). Recall that a bounded net $(e_\alpha)_\alpha$ in $\mathfrak{A}$ is called a bounded left (right) approximate identity for $X,$ if $\pi_\ell(e_\alpha, x)\rightarrow x\;(\pi_r(x,e_\alpha) \rightarrow x),$ for every $x\in X.$ Furthermore, if $\mathcal{A}$ has the unit element $e$ and
$$\pi_\ell(e,x)=\pi_r(x,e)=x \qquad (x\in X),$$
we say that $X$ is a unital $\mathcal{A}$-bimodule.

Let $X$ be a Banach $\mathcal{A}$-bimodule.
  A \emph{derivation} from $\mathcal{A}$ into $X$ is a bounded linear mapping $D:\mathcal{A}\rightarrow X$ such that $$D(ab)=aD(b)+D(a)b \text{ for all } a,b\in \mathcal{A}.$$
The space of continuous derivations from $\mathcal{A}$ into $X$ is denoted by $Z^1(\mathcal{A},X)$.
Easy examples of derivations are the inner derivations, which are given for each $x\in X$ by
$$\delta_x(a)=ax-xa \text{ for all }a\in \mathcal{A}.$$
The space of inner derivations from $\mathcal{A}$ into $X$ is denoted by $N^1(\mathcal{A},X)$.
The Banach algebra $\mathcal{A}$ is said to be \emph{amenable}, when for every Banach $\mathcal{A}$-bimodule $X$, the inner derivations are only derivations existing from $\mathcal{A}$ into $X^*$. It is clear that $\mathcal{A}$ is amenable if and only if $H^1(\mathcal{A},X^*)=Z^1(\mathcal{A},X^*)/ N^1(\mathcal{A},X^*)=\{0\},$ for every Banach $\mathcal{A}$-bimodule $X.$
The concept of weak amenability was first introduced by Bade, Curtis and Dales in [2] for commutative Banach algebras, and was extended to the noncommutative case by Johnson, see [13].
For Banach $A-bimodule$ $X$, the quotient space $H^1(A,X)$ of all continuous derivations from $A$ into $X$ modulo the subspace of inner derivations is called the first cohomology group of $A$ with coefficients in $X$.\\
A Banach algebra $\mathcal{A}$ is said to be \emph{weakly amenable}, if every derivation from $\mathcal{A}$ into $\mathcal{A}^*$ is inner. Similarly, $\mathcal{A}$ is weakly amenable if and only if $H^1(\mathcal{A},\mathcal{A}^*)=Z^1(\mathcal{A},\mathcal{A}^*)/ N^1(\mathcal{A},\mathcal{A}^*)=\{0\}$.
\vspace{0.4cm}
\section{Arens regularity and weak amenability}
\vspace{0.4cm}
\begin{theorem}\label{thm1}
Let $X$ be a Banach $\mathcal{A}$-bimodule and $X^\ast$ factors on the left with respect to $\mathcal{A}.$ If $\mathcal{A}\mathcal{A}^{\ast\ast}\subseteq Z_1(\pi_\ell),$ then $Z_1 (\pi_\ell)=\mathcal{A}^{\ast\ast}.$
\end{theorem}

\begin{proof}
Let $a'' \in \mathcal{A}^{\ast\ast}$ and $x'' \in X^{\ast\ast}.$ Also suppose that $(x''_\alpha)_\alpha$ is a convergent net to  $x''$ in $\sigma(X^{\ast\ast},X^{\ast}).$ Since $X^*$ factors on the left, for every $x' \in X^{\ast},$ there are $a \in \mathcal{A}$ and $y' \in X^{\ast}$ such that $x' = y'a.$ Since $aa'' \in Z_1(\pi_\ell),$ we have
\begin{align*}
  \lim_\alpha \langle \pi_{\ell}^{\ast\ast\ast} (a'',x''_\alpha), x' \rangle &= \lim_\alpha \langle \pi_{\ell}^{\ast\ast\ast} (a'',x''_\alpha), y'a \rangle = \lim_\alpha \langle \pi_{\ell}^{\ast\ast\ast} (aa'',x''_\alpha), y' \rangle\\ &= \langle \pi_{\ell}^{\ast\ast\ast} (aa'',x''), y' \rangle=\langle \pi_{\ell}^{\ast\ast\ast} (a'',x''), x' \rangle.
\end{align*}
It follows that $\pi_{\ell}^{\ast\ast\ast} (a'',x''_\alpha) \rightarrow \pi_{\ell}^{\ast\ast\ast} (a'',x'')$ in $\sigma(X^{\ast\ast},X^{\ast}),$ and so $a'' \in Z_1(\pi_\ell)$. Thus $Z_1(\pi_\ell)=\mathcal{A}^{\ast\ast}.$
\end{proof}
\vspace{0.4cm}
\begin{theorem}
Let $X$ be a Banach $\mathcal{A}$-bimodule and $X^\ast$ factors on the left with respect to $\mathcal{A}.$ If $\mathcal{A}X^{\ast\ast}\subseteq Z_1(\pi_r),$ then $Z_1(\pi_r)=X^{\ast\ast}.$
\end{theorem}

\begin{proof}
The proof is similar to the Theorem \ref{thm1}.
\end{proof}
\vspace{0.4cm}
\begin{corollary}
Let $\mathcal{A}$ be a Banach algebra and $\mathcal{A}^\ast$ factors on the left. Then $\mathcal{AA}^{\ast\ast}\subseteq Z_1 (\mathcal{A}^{\ast\ast})$ if and only if $\mathcal{A}$ is Arens regular.
\end{corollary}
\vspace{0.4cm}
By \cite{15}, for a locally compact group $G,$ its group algebra $L^1(G)$ factors on the right \ff $L^1(G)$ factors on the left \ff $G$ is discrete. On the other hand, $L^1(G)$ is Arens regular \ff $G$ is finite. Consequently, the following result is immediate.
\vspace{0.4cm}
\begin{corollary}
Let $G$ be any group. Then we have $\ell^1(G)\ell^\infty (G)^\ast \subseteq \ell^1 (G)$ \ff $G$ is finite.
\end{corollary}

\vspace{0.4cm}

\begin{theorem}
Let $X$ be a Banach $\mathcal{A}$-bimodule and $\mathcal{A}^{\ast\ast}$ has a \textnormal{BAI} for $X^{\ast\ast}.$ Then $X^{\ast\ast}$ is a unital Banach $\mathcal{A}^{\ast\ast}$-bimodule.
\end{theorem}

\begin{proof}
Let $(e''_\alpha)_\alpha$ be a \textnormal{BAI} for $X^{\ast\ast}.$ By passing to a subnet, we may assume that $e''_\alpha \rightarrow e''$ in $\sigma(\mathcal{A}^{\ast\ast},\mathcal{A}^{\ast}),$ for some $e''$ in $\mathcal{A}^{\ast\ast}.$ Then for every $x''\in X^{\ast\ast}$ and $x' \in X^{\ast},$ we have
$$\langle \pi_{\ell}^{\ast\ast\ast} (e'',x''), x' \rangle = \lim_\alpha \langle \pi_{\ell}^{\ast\ast\ast} (e_\alpha'',x''), x' \rangle = \langle x'', x'\rangle.$$
Thus $\pi_{\ell}^{\ast\ast\ast} (e'',x'') = x''.$

Now for every $x \in X$ we have
$$\pi_r^{\ast\ast\ast} (x,e'') = \text{weak}^\ast-\lim_\alpha \pi_r^{\ast\ast\ast} (x,e''_\alpha) =x. $$
Take $x'' \in X^{\ast\ast}$ arbitrary and $(x_\beta)_\beta \subseteq X$ such that $x_\beta \rightarrow x''$ in $\sigma(X^{\ast\ast},X^{\ast}).$ Then
$$\pi_r^{\ast\ast\ast} (x'',e'') = \text{weak}^\ast-\lim_\beta \pi_r^{\ast\ast\ast} (x_\beta,e'') = \text{weak}^\ast-\lim_\beta x_\beta = x''. $$
Consequently $\pi_r^{\ast\ast\ast} (x'',e'') = x''.$
\end{proof}
\vspace{0.4cm}
\begin{corollary}
Let $\mathcal{A}$ be a Banach algebra such that $\mathcal{A}^{\ast\ast}$ is amenable. Then $\mathcal{A}^{\ast\ast}$ has the identity element.
\end{corollary}
\vspace{0.4cm}
\begin{proof}
The proof is clear by the preceding theorem and the fact that every amenable Banach algebra has a BAI, see {[9, Proposition 1.6]}.
\end{proof}
\vspace{0.4cm}
\begin{theorem}
Let $X$ be a left Banach $\mathcal{A}$-module. Then we have the following assertions:
\begin{enumerate}
\item[(i)] $\pi_\ell^{\ast\ast\ast\ast} (X^{\ast},\mathcal{A}^{\ast\ast}) \subseteq X^{\ast}$ \ff $Z_1(\pi_\ell) = \mathcal{A}^{\ast\ast}.$
\item[(ii)] The mapping $x'' \mapsto \pi_\ell^{\ast\ast}(x'',x')$ from $X^{\ast\ast}$ into $\mathcal{A}^{\ast}$ is weak$^{\ast}$ to weak continuous for all $x' \in X^\ast$ \ff $Z_1(\pi_\ell) = \mathcal{A}^{\ast\ast}.$
\item[(iii)] If for every $x''' \in X^{\ast\ast\ast},$ the mapping $x'' \mapsto \pi_\ell^{\ast\ast\ast\ast\ast}(x'',x''')$ from $X^{\ast\ast}$ into $\mathcal{A}^{\ast\ast\ast}$ is weak$^{\ast}$ to weak$^{\ast}$ continuous, then the mapping $x'' \mapsto \pi_\ell^{\ast\ast\ast}(a'',x'')$ from $X^{\ast\ast}$ into $X^{\ast\ast}$ is weak$^{\ast}$ to weak continuous, and so $Z_1(\pi_\ell) = \mathcal{A}^{\ast\ast}.$
\end{enumerate}
\end{theorem}

\begin{proof}
(i): Suppose that $a'' \in \mathcal{A}^{\ast\ast}$ and $(x''_\alpha)_\alpha\subseteq X^{\ast\ast},$ such that $x''_\beta \rightarrow x''$ in the weak$^*$-topology. Let $x' \in X^\ast$ and $x \in X.$ Since $\pi_\ell^{\ast\ast\ast\ast} (X^{\ast},\mathcal{A}^{\ast\ast}) \subseteq X^{\ast}$, $\pi_\ell^{\ast\ast\ast\ast} (x',a'') \in X^{\ast}$.
Thus we have
\begin{align*}
  \lim_\alpha \langle \pi_{\ell}^{\ast\ast\ast} (a'',x_\alpha''), x' \rangle &= \lim_\alpha \langle a'',\pi_{\ell}^{\ast\ast} (x''_\alpha,x')\rangle = \lim_\alpha \langle \pi_{\ell}^{\ast\ast\ast\ast\ast} (x''_\alpha,x'), a'' \rangle\\ &=\lim_\alpha \langle x''_\alpha,\pi_{\ell}^{\ast\ast\ast\ast} (x',a'')\rangle=
  \langle x'', \pi_{\ell}^{\ast\ast\ast\ast} (x',a'')\rangle \\&= \langle \pi_{\ell}^{\ast\ast\ast} (a'',x'') , x'\rangle.
\end{align*}
žConsequently, we have  $\pi_{\ell}^{\ast\ast\ast} (a'',x_\alpha'')\rightarrow \pi_{\ell}^{\ast\ast\ast} (a'',x'')$ in weak$^*$ topology in $X^{**}$, and so $a''\in Z_1(\pi_\ell)$.\\
Conversely, suppose that $Z_1(\pi_\ell) = \mathcal{A}^{\ast\ast}.$ Take $(x''_\alpha)_\alpha \subseteq X^{\ast\ast}$ such that $x''_\beta \rightarrow x''$ in $(X^{**},\sigma(X^{\ast\ast},X^{\ast})).$ Then
\begin{align*}
  \lim_\alpha \langle \pi_{\ell}^{\ast\ast\ast\ast} (x',a''),x''_\alpha\rangle &= \lim_\alpha \langle \pi_{\ell}^{\ast\ast\ast} (a'',x''_\alpha),x'\rangle\\& = \langle \pi_{\ell}^{\ast\ast\ast} (a'',x''),x'\rangle\\& = \langle \pi_{\ell}^{\ast\ast\ast\ast} (x',a''),x''_\alpha\rangle.
\end{align*}
It follows that $\pi_{\ell}^{\ast\ast\ast\ast} (x',a'') \in (X^{\ast\ast},\text{weak}^\ast)^\ast = X^\ast.$ Consequently, we have
$$\pi_{\ell}^{\ast\ast\ast\ast} (X^\ast,A^{\ast\ast}) \subseteq X^{\ast}.$$

(ii): Let the mapping $x''\mapsto \pi_\ell^{\ast\ast}(x'',x')$ be weak$^\ast$ to weak continuous. Suppose that $a'' \in \mathcal{A}^{\ast\ast}$ and $(x''_\alpha)_\alpha \subseteq X^{\ast\ast},$ such that $x''_\beta \rightarrow x''$ in $(X^{**},\sigma(X^{\ast\ast},X^{\ast})).$ Then for every $x' \in B^\ast,$ we have
\begin{align*}
  \lim_\alpha \langle \pi_{\ell}^{\ast\ast\ast} (a'',x''_\alpha), x' \rangle &= \lim_\alpha \langle a'',  \pi_{\ell}^{\ast\ast} (x''_\alpha, x' )\rangle \\&=
  \langle a'',  \pi_{\ell}^{\ast\ast} (x'', x' )\rangle \\&= \langle \pi_{\ell}^{\ast\ast\ast} (a'',x''), x' \rangle.
\end{align*}

(iii): Let $a'' \in \mathcal{A}^{\ast\ast}$ and $(x''_\alpha)_\alpha \subseteq X^{\ast\ast},$ such that $x''_\beta \rightarrow x''$ in $(X^{**},\sigma(X^{\ast\ast},X^{\ast})).$ Then for every $x'''\in X^{\ast\ast\ast}$ we have
\begin{align*}
  \lim_\alpha \langle x''', \pi_{\ell}^{\ast\ast\ast} (a'',x''_\alpha)\rangle &= \lim_\alpha \langle \pi_{\ell}^{\ast\ast\ast\ast} (x''',a''), x''_\alpha \rangle = \lim_\alpha\langle x''_\alpha , \pi_{\ell}^{\ast\ast\ast\ast} (x''', a'' )\rangle \\&= \lim_\alpha \langle \pi_{\ell}^{\ast\ast\ast\ast} (x''_\alpha,x'''), a'' \rangle = \langle \pi_{\ell}^{\ast\ast\ast\ast\ast} (x'',x'''), a'' \rangle\\& = \langle x''' , \pi_{\ell}^{\ast} (a'',x'') \rangle.
\end{align*}
It follows that $\pi_{\ell}^{\ast\ast\ast} (a'',x''_\alpha) \rightarrow \pi_{\ell}^{\ast\ast\ast} (a'',x'')$ weakly as required.
\end{proof}
\vspace{0.4cm}
\begin{corollary}
By each of the following conditions, a Banach algebra $\mathcal{A}$ is Arens regular.
\begin{enumerate}
\item[(i)] $\mathcal{A}$ is a left ideal in its second dual and $\mathcal{A}^*$ factors on the left.
\item[(ii)] The mapping $a'' \mapsto a''a'$ is $(\mathcal{A}^{**}, \sigma(\mathcal{A}^{**},\mathcal{A}^{*}))$ to $(\mathcal{A}^{*}, \sigma(\mathcal{A}^{*},\mathcal{A}))$ continuous, for every $a' \in \mathcal{A}^{\ast}.$
\item[(iii)] The mapping $a'' \mapsto a''a'''$ is $(\mathcal{A}^{**}, \sigma(\mathcal{A}^{**},\mathcal{A}^{*}))$ to $(\mathcal{A}^{***}, \sigma(\mathcal{A}^{***},\mathcal{A}^{**}))$ continuous, for every $a''' \in \mathcal{A}^{***}.$
\end{enumerate}
\end{corollary}
\vspace{0.4cm}
The next result provides some conditions under which weak amenability is inherited by a Banach algebra $\mathcal{A}$ from the second dual $\mathcal{A}^{**}.$
\vspace{0.4cm}
\begin{definition}
Let $X$ be a Banach $\mathcal{A}$-bimodule and $x'\in X^*.$ We say that $x'$ has RW$^*$WC-\emph{property} with respect to $\mathcal{A}'',$ if for every $(a''_\alpha)_\alpha \subseteq \mathcal{A}^{**}$ such that $a''_\alpha x' \rightarrow 0$ in weak$^*$ topology, it follows that $a''_\alpha x' \rightarrow 0$ weakly.
\end{definition}
\vspace{0.4cm}
\begin{theorem}
For a Banach $\mathcal{A}$-bimodule $X,$ the following assertions are equivalent:
\begin{enumerate}
\item[(i)] Every $x' \in X^*$ has the RW$^*$WC-property with respect to $\mathcal{A}''$.
\item[(ii)] $Z_1(\pi_r)=X^{**}.$
\end{enumerate}
\end{theorem}
\vspace{0.4cm}
\begin{proof}
(i)$\rightarrow$(ii): Let $(a''_\alpha)_\alpha \subseteq \mathcal{A}^{**}$ such that $a''_\alpha \rightarrow a''$ in $(\mathcal{A}^{**},\sigma(\mathcal{A}^{**},\mathcal{A}^{*})).$ Then for every $x \in X$ and $x' \in X^*$ we have
$$\lim_\alpha \langle a''_\alpha x',x \rangle = \lim_\alpha \langle a''_\alpha , x'x \rangle= \langle a'' , x'x\rangle= \langle a''x', x\rangle.$$
It follows that $a''_\alpha x' \rightarrow a''x'$ in weak$^*$ topology. Since every $x' \in X^*$ has the RW$^*$WC-property with respect to $\mathcal{A}^{**}$, for every $x''\in X^{**}$, we have

$\lim_\alpha \langle x''a''_\alpha ,x' \rangle = \lim_\alpha \langle x'' ,a''_\alpha x' \rangle= \langle x'' , a''x'\rangle= \langle x''a'', x'\rangle.$
It follows that $a''_\alpha x''\rightarrow a'' x''$ in weak$^*$ topology. Hence $Z_1(\pi_r)=X^{**}.$\\

(ii)$\rightarrow$(i): Let $(a''_\alpha)_\alpha \subseteq \mathcal{A}^{**}$ such that $a''_\alpha \rightarrow 0$ in $(\mathcal{A}^{**},\sigma(\mathcal{A}^{**},\mathcal{A}^{*}))$ and $a''_\alpha x' \rightarrow 0$ in weak$^*$ topology where $x^\prime\in X^*$. Since $Z_1(\pi_r)=X^{**},$ for every $x^{\prime\prime}\in X^{**}$, we have
$$\langle x^{\prime\prime},a''_\alpha x' \rangle=\langle x^{\prime\prime}a''_\alpha , x' \rangle=0$$
It follows that $a''_\alpha x'\rightarrow 0$ in weak topology, and so   $x' \in X^*$ has the RW$^*$WC-property with respect to $\mathcal{A}''$.
\end{proof}
\vspace{0.4cm}

\begin{theorem}
Let $\mathcal{A}$ be a Banach algebra and the mapping $a'' \mapsto a''a'''$ from $\mathcal{A}^{**}$ into $\mathcal{A}^{***}$ is $\sigma(\mathcal{A}^{**},\mathcal{A}^{*})$ to $\sigma(\mathcal{A}^{***},\mathcal{A}^{**})$ continuous, for every $a''' \in A^{\ast\ast\ast}.$ Then weak amenability of $\mathcal{A}^{\ast\ast}$ implies this property for $\mathcal{A}.$
\end{theorem}

\begin{proof}
Assume that $D:\,\mathcal{A} \longrightarrow \mathcal{A}^\ast$ is a derivation. We show that $D'':\,\mathcal{A}^{\ast\ast} \longrightarrow \mathcal{A}^{\ast\ast\ast}$ is also a derivation. Take $a''$ and $b''$ in $\mathcal{A}^{\ast\ast}$ arbitrary and then take $(a_\alpha)_\alpha$ and $(b_\beta)_\beta$ two nets in $\mathcal{A}$ convergent to $a''$ and $b'',$ respectively, in $\sigma(\mathcal{A}^{\ast\ast},\mathcal{A}^{\ast})$-topology. Since by assumption the mapping $a'' \mapsto a''a'''$ is weak$^{\ast}$ to weak$^{\ast}$ continuous, then we have
$$\lim_\alpha \lim_\beta a_\alpha D(b_\beta)= a''D''(b'').$$
Thus $D''(a''b'') = \lim_\alpha \lim_\beta D(a_\alpha b_\beta)= \lim_\alpha \lim_\beta [a_\alpha D(b_\beta)+D(a_\alpha) b_\beta] = a''D''(b'') + D(a'') b''.$
Since $\mathcal{A}^{\ast\ast}$ is weakly amenable, there exists $a'''\in \mathcal{A}^{***}$ such that $D''=\delta_{a'''}.$ Consequently $D=\delta_{a'},$ where $a'$ is the restriction of $a'''$ on $\mathcal{A}.$
\end{proof}
\vspace{0.4cm}
Let $A$ be a Banach Algebra. A dual Banach $A-bimodule$ $X$ is called normal if, for each $x\in X$ the map $a\rightarrow ax$ and $a\rightarrow xa$ from $A$ into $X$ is $weak^*-to-weak^*$ continuous.\\
A dual Banach algebra $A$ is Connes-amenable if, for every normal, dual Banach $A-bimodule$ $X$, every $weak^*-to-weak^*$ continuous derivation $D\in Z^1(A,X)$ is inner. Then we write $H^1_{w^*}(A,X)=0$.\\
Let $A$ be a Banach Algebra. We denote by $Z_{w^*}^1(A,X)$ for every $weak^*-to-weak^*$ continuous derivation $D\in Z^1(A,X)$ where $X$ is a Banach $A-bimodule$. Thus we can write first $weak^*$ cohomological group on $A$ as follows
$$H_{w^*}^1(\mathcal{A},X^*)=Z_{w^*}^1(\mathcal{A},X^*)/ N_{w^*}^1(\mathcal{A},X^*).$$

\vspace{0.4cm}
\begin{theorem}
Let $\mathcal{A}$ be an Arens regular Banach algebra. If  $A$ is weakly amenable, it follows that $H^1(\mathcal{A}^{\ast\ast}, \mathcal{A}^\ast)=0$.
\end{theorem}

\begin{proof} Let  $D\in Z^1_{w^*}(A^{**},A^*)$. Take $D_1=D\mid_A$. Then it is clear that $D_1\in Z^1(A,A^*)$. Since
$H^1(\mathcal{A}, \mathcal{A}^\ast)=0$, there is $a^\prime\in A^*$ such that $D_1=\delta_{a^\prime}$ in $A$.  Assume that  $a'' \in \mathcal{A}^{\ast\ast}$ and $(a_\alpha)_\alpha\subseteq A^{\ast\ast},$ such that $a_\alpha\rightarrow a''$ in the weak$^*$-topology. Then Since $A$ is Arens regular, by using Theorem 2.7 part (i), for every $a^\prime\in A^*$, we have $a^\prime a^{\prime\prime}\in A^*$, and so $w^*-\lim_\alpha (a^\prime a_\alpha)=a^\prime a^{\prime\prime}$. Then we have
$$D(a^{\prime\prime})=D(w^*-\lim_\alpha a_\alpha)=w^*-\lim_\alpha D(a_\alpha)=w^*-\lim_\alpha \delta_{a^\prime}(a_\alpha)=w^*-\lim_\alpha (a^\prime a_\alpha-a_\alpha a^\prime)$$$$=a^\prime a^{\prime\prime}-a^{\prime\prime}a^\prime=\Delta_{a^\prime}(a^{\prime\prime}).$$
Thus $D=\Delta_{a^\prime}$ on $A^{**}$, and so $H^1(\mathcal{A}^{\ast\ast}, \mathcal{A}^\ast)=0$.
\end{proof}
\vspace{0.4cm}

\begin{example}
 Consider the algebra $c_0=(c_0,.)$ is the collection of all sequences of scalars that convergence to $0$, with the some vector space operations and norm as $\ell^\infty$. By using [5, Example 2.6.22 ] we know that $c_0=(c_0,.)$ is Arens regular, and also  by [6], $c_0$ is weakly amenable. Thus by using preceding  theorem , it is clear that  $H^1(c_0^{\ast\ast}, c_0^\ast)=0$.
\end{example}
\vspace{0.4cm}

\noindent {\bf Problem.}  Let $\mathcal{A}$ be a Banach algebra and the mapping $a'' \mapsto a''a'$ from $\mathcal{A}^{\ast\ast}$ into $\mathcal{A}^\ast$ is weak$^\ast$ to weak continuous, for every $a' \in \mathcal{A}^{\ast}.$ Then,   $H^1(\mathcal{A}, \mathcal{A}^\ast)=0$ \ff $H^1(\mathcal{A}^{\ast\ast}, \mathcal{A}^\ast)=0$.

\vspace{0.4cm}

\bibliographystyle{amsplain}

\begin{thebibliography}{99}
\bibitem{1} R. E.  Arens, {\it The adjoint of a bilinear operation}, Proc. Amer. Math. Soc. {\bf 2} (1951), 839-848.
\bibitem{3} W.  G. Bade, P.C. Curtis and H.G. Dales, {\it Amenability and weak amenability for Beurling and Lipschitz algebra}, Proc. Lodon Math. Soc. {\bf 50} (1987) 359-377.

\bibitem{3} J. Baker, A.T. Lau, J.S. Pym {\it Module homomorphism and topological centers associated with
 weakly sequentially compact Banach algebras}, Journal of Functional Analysis. {\bf 158} (1998), 186-208.
\bibitem{4} F. F. Bonsall, J. Duncan, {\it Complete normed algebras}, Springer-Verlag, Berlin 1973.

\bibitem{5} H. G. Dales, {\it Banach algebra and automatic continuity}, Oxford 2000.
\bibitem{6} H. G. Dales, F. Ghahramani, N. Gr{\o}nb{\ae}k {\it Derivation into iterated duals of Banach algebras} Studia Math. {\bf 128} 1 (1998), 19-53.
\bibitem{7} H. G. Dales, A. Rodrigues-Palacios, M.V. Velasco, {\it The second transpose of a derivation}, J. London. Math. Soc. {\bf2} 64 (2001) 707-721.


\bibitem{8} N. Dunford, J. T. Schwartz, {\it Linear operators.I},
Wiley, New york 1958.
\bibitem{8} H. Farhadi, F.Ghahramani,  {\it Involutions on the second duals of group algebras and a multiplier problem}, Proceedings of the Edinburgh Mathematical Society  {\bf 50}  (2007), 153-161.
\bibitem{10} M. Eshaghi Gordji, M. Filali, {\it Weak amenability of the second dual of a Banach algebra}, Studia Math. {\bf182} 3 (2007), 205-213.
\bibitem{11}  F. Gourdeau, {\it Amenability of Lipschits algebra}, Math. Proc. Cambridge. Philos. Soc. {\bf112} (1992), 581-588.
\bibitem{12} E. Hewitt, K. A. Ross,  {\it Abstract harmonic analysis}, Springer, Berlin, Vol I 1963.

\bibitem{10} B. E. Johoson,  {\it Derivation from $L^1(G)$ into $L^1(G)$ and $L^\infty(G)$}, Harmonic analysis. Luxembourg 1987, 191-198 Lecture Note in Math., {\bf 1359}, Springer, Berlin, 1988. MR  {\bf 90a}:46122.

\bibitem{10} B. E. Johoson,  {\it Weak amenability of group algebra}, Bull. Lodon. Math. Soc. {\bf 23}(1991), 281-284.
\bibitem{14}  A. T. Lau, V. Losert, {\it On the second Conjugate Algebra of locally
compact groups}, J. London Math. Soc.  {\bf 37} (2)(1988),
464-480.
\bibitem{15} A. T. Lau, A. \"{U}lger, {\it Topological center of certain dual
algebras}, Trans. Amer.  Math. Soc. {\bf 348} (1996), 1191-1212.

\bibitem{16} S. Mohamadzadih, H. R. E. Vishki, {\it Arens regularity of module actions and the second adjoint of a derivation}, Bulletin of the Australian Mathematical Society {\bf77} (2008), 465-476.
\bibitem{17}  J. S. Pym, {\it The convolution of functionals on spaces of bounded functions},
 Proc. London Math Soc.  {\bf 15} (1965), 84-104.

\bibitem{18} A. \"{U}lger, {\it Some stability properties of Arens regular bilinear operators}, Proc. Amer. Math. Soc. (1991) {\bf 34}, 443-454.
\bibitem{19}  A. \"{U}lger, {\it Arens regularity of weakly sequentialy compact Banach algebras},
Proc. Amer. Math. Soc. {\bf 127} (11) (1999), 3221-3227.
\bibitem{20} P. K. Wong, {\it The second conjugate algebras of
Banach algebras}, J. Math. Sci. {\bf 17} (1) (1994), 15-18.
\bibitem{21} Y. Zhing,  {\it Weak amenability of module extentions of Banach algebras}, Trans. Amer. Math. Soc. {\bf 354} (10) (2002), 4131-4151.
\bibitem{22} Y. Zhing,  {\it Weak amenability of a class of Banach algebra}, Cand. Math. Bull. {\bf 44} (4) (2001) 504-508.\\
\end{thebibliography}

\end{document}